\documentclass[12pt,a4paper,oneside]{amsart}

\usepackage{amsfonts, amsmath, amssymb, amsthm, hyperref}
\usepackage{anysize}

\newtheorem{thm}{Theorem}
\newtheorem*{thm*}{Theorem}

\theoremstyle{definition}

\theoremstyle{remark}
\newtheorem*{rem}{Remark}

\renewcommand{\int}{\mathop{\rm int}}

\renewcommand{\epsilon}{\varepsilon}

\begin{document}

\title{A simpler proof of the Boros--F\"uredi--B\'ar\'any--Pach--Gromov theorem}

\author{Roman~Karasev}
\thanks{This research is supported by the Dynasty Foundation, the President's of Russian Federation grant MK-113.2010.1, the Russian Foundation for Basic Research grants 10-01-00096 and 10-01-00139, the Federal Program ``Scientific and scientific-pedagogical staff of innovative Russia'' 2009--2013}
\email{r\_n\_karasev@mail.ru}
\address{ Roman Karasev, Dept. of Mathematics, Moscow Institute of Physics and Technology, Institutskiy per. 9, Dolgoprudny, Russia 141700}


\begin{abstract}
A short and almost elementary proof of the Boros--F\"uredi--B\'ar\'any--Pach--Gromov theorem on the multiplicity of covering by simplices in $\mathbb R^d$ is given.
\end{abstract}

\maketitle

Let us give a proof of the Boros--F\"uredi--B\'ar\'any--Pach--Gromov theorem~\cite{ba1982,bf1984,pach1998,grom2010} that is actually the ``decoded'' and refined proof from~\cite{grom2010} (see also~\cite[Section~2]{foxetal2010} for a similar proof in the two-dimensional case). Unlike the proof in~\cite{grom2010}, the only topological notion that is used here is the degree of a piece-wise smooth map. 

Consider a set of $d+1$ absolutely continuous
probability measures $\mu_0,\mu_1,\ldots, \mu_d$ on $\mathbb R^d$. Define a \emph{random simplex} of dimension $k$ as a simplex spanned by $k+1$ points $x_{d-k},\ldots, x_d \in \mathbb R^d$, where the point $x_i$ is distributed according to the measure $\mu_i$. The following theorem estimates from below the maximum \emph{simplicial depth}~\cite{liu1990} over all points in $\mathbb R^d$, i.e. the probability for this point to be covered by a random $d$-simplex:

\begin{thm}
\label{bfbcg}
Under the above assumptions there exists a point $c\in \mathbb R^d$ such that the probability for a random $d$-simplex to contain $c$ is
$$
\ge p_d=\frac{1}{(d+1)!}.
$$
\end{thm}

Note that in~\cite{grom2010} a stronger result is proved: the maps $\Delta^N\to Y$ of a simplex with measure to a smooth manifold were considered. Here we give the statement of Theorem~\ref{bfbcg} that is closer to the original theorems in~\cite{ba1982,bf1984,pach1998}.

\begin{proof}[Proof of Theorem~\ref{bfbcg}]
We assume that $\mathbb R^d$ is contained in its one-point compactification $S^d = \mathbb R^d\cup\{\infty\}$.

Assume the contrary. Take some small $\varepsilon>0$. Consider a fine enough finite triangulation $Y$ of $S^d$ with one vertex at $\infty$ so that for any $0<k\le d$ and any $k$-face $\sigma$ of $Y$ the probability of a random $(d-k)$-simplex $x_kx_{k+1}\ldots x_d$ to intersect $\sigma$ is $<\varepsilon$. Here and below we always assume that $\mu_i$ is the distribution of $x_i$. 
To make such a triangulation it is sufficient to take a large enough ball $B$ so that at least $1-\varepsilon$ of every measure is inside $B$. Then we take the simplices of $Y$ that intersect $B$ small enough, other simplices may be arbitrary. From the absolute continuity it follows that for small enough simplices the probabilities become arbitrarily small, and for simplices in $\mathbb R^d\setminus B$ they are $<\varepsilon$ by the choice of $B$. 

Consider a $(d+1)$-dimensional simplicial complex $Y*0$ (the cone over $Y$ with apex $0$). Now we are going to build a (piece-wise smooth) map $f : (Y*0)^{(d)}\to S^d$ (from the $d$-skeleton) which is ``economical'' with respect to the measures $\mu_i$ (this phrase will be clarified below), and coincides with the identification $Y=S^d$ on $Y\subset (Y*0)^{(d)}$.

Proceed by induction:

\begin{itemize}
\item 
Map $0$ to $\infty\in S^d$;
\item 
For any vertex $v\in Y$ map $[v0]$ to an open ray starting from $v$ (and ending at $\infty\in S^d$) so that the probability for a random $(d-1)$-simplex $x_1\ldots x_d$ to meet $f([v0])$ is $< p_d$. This is possible because a simplex $x_0x_1\ldots x_d$ contains $v$ iff the $(d-1)$-simplex $x_1\ldots x_d$ intersects the ray from $v$ opposite to $x_0-v$. Since the probability for a random $d$-simplex to contain $v$ is $<p_d$, for some of such rays the corresponding probability is also $<p_d$.
\item 
Step to the $k$-skeleton of $Y*0$ as follows. Let $\sigma=v_1\ldots v_k0$ be a $k$-simplex of $Y*0$. The map $f$ is already defined for $\partial\sigma$. We know that the probability for a random $(d-k+1)$-simplex $x_{k-1}\ldots x_d$ to meet some $f(v_1\ldots\hat{v_i}\ldots v_k0)$ ($i=1,\ldots,k$) is $<(k-1)!p_d$, and the probability to meet $f(v_1\ldots v_k)$ is $<\varepsilon$. If $\varepsilon$ is chosen small enough we see that a random $(d-k+1)$-simplex $x_{k-1}\ldots x_d$ intersects $f(\partial\sigma)$ with probability $<k!p_d$. 

There exist a point $x_{k-1}$ not in $f(\partial\sigma)$ such that the probability for $x_{k-1}x_k\ldots x_d$ (with random last $d-k+1$ points) to meet $f(\partial\sigma)$ is $<k!p_d$; here the independence of the distributions of vertices is essential. Let us define the map $f$ on the simplex $\sigma$ treated as a join $\partial\sigma*c$ so that $c$ is mapped to $\infty\in S^d$, and every segment $[vc]$ ($v\in\partial\sigma$) is mapped to the infinite ray from $f(v)$ in the direction opposite to $x_{k-1}-v$. More explicitly: map $[vc]$ to $[f(v), x_{k-1}]$ first; then apply the inversion with center $x_{k-1}$ and radius $|x_{k-1} - f(v)|$ that maps $[f(v), x_{k-1}]$ to $[f(v), \infty]$; if $f(v)=\infty$ then map $[vc]$ to the point $\infty$. Now the probability for a random $(d-k)$-simplex to intersect $f(\sigma)$ is $<k!p_d$. 
\end{itemize}

Finally for any $d$-simplex $\sigma$ of $Y$ we have that the boundary of the cone $\sigma*0$ is mapped so that 
$$
\mu_d(f(\partial(\sigma*0)))<(d+1)!p_d=1,
$$ 
if we again use small enough $\varepsilon$. Therefore $f(\partial(\sigma*0))\neq S^d$ and the restriction $f|_{\partial(\sigma*0)}$ has zero degree. By summing up the degrees (the $d$-faces of $(\partial\sigma)*0$ go pairwise and cancel, because $Y$ is a triangulation) we see that the map $f|Y$ has even degree  but it is the identity map, which is a contradiction.  
\end{proof}

This theorem can be sharpened (following~\cite{grom2010}) if two of the measures coincide.

\begin{thm}
\label{bfbcg2}
If some two measures coincide then the bound in Theorem~\ref{bfbcg} can be improved to
$$
p'_d=\frac{2d}{(d+1)!(d+1)}.
$$
\end{thm}

\begin{proof}
Assume $\mu_{d-1} = \mu_d$. We proceed in the same way building $f : (Y*0)^{(d)}\to S^d$, but we slightly change the construction on the last step.

On the last step we have a $(d-1)$-simplex $\sigma$ of $Y$, and $f$ is already defined for $\partial (\sigma*0)$ so that the probability for a random segment $[x_{d-1}x_d]$ to intersect $D=f(\partial(\sigma*0))$ is $<d!p'_d = \frac{2d}{(d+1)^2}$.

We are going to extend $f$ to $\sigma*0$ so that its image $f(\sigma*0)$ mod $2$  has measure $<\frac{1}{d+1}$ (we use the measure $\mu_{d-1}=\mu_d$). Here the image mod $2$ is the set of points in $\mathbb R^d$ that are covered by $f(\sigma*0)$ odd number of times.

It can be easily seen that $D$ ``partitions'' $\mathbb R^d$ into two parts $A$ and $B$ characterized by the following property: any generic piece-wise linear path from $A$ to $B$ meets $D$ odd number of times, and any generic piece-wise linear path with both ends in $A$ (or both in $B$) meets $D$ even number of times. The sets $A$ and $B$ are the only possibilities of image of $f(\sigma *0)$ mod $2$, because the covering parity of $f|_{\sigma*0}$ changes only at crossing with $f(\partial (\sigma*0))= D$.

If $\mu_d(A)=x$ and $\mu_d(B)=1-x$ then the probability for a random segment $[x_{d-1}x_d]$ (recall that $\mu_{d-1}=\mu_d$) to meet $D$ is at least $2x(1-x)$, that is
$$
\frac{2d}{(d+1)^2} > 2x(1-x). 
$$
It follows easily that in this case either $x$ or $1-x$ is $<\frac{1}{d+1}$ and we can define $f$ as required again.

Now let $\tau$ be a $d$-simplex of $Y$. For any its facet $\sigma\subset \tau$ the image $f(\sigma*0)$ mod $2$ has measure $<\frac{1}{d+1}$. It follows that $f(\partial\tau * 0)$ mod $2$ has measure $<1$; because if a point has odd number of preimages in $\partial \tau * 0$ then it has odd number of preimages in some $\sigma * 0$ (this is true if the preimage is not in the $(d-1)$-skeleton, but such exceptions correspond to a zero measure set). Since the measure of $\tau$ itself if $<\varepsilon$, we see that $f(\partial(\tau*0))$ mod $2$ has measure $<1$ and therefore $f$ has even degree on $\partial(\tau*0)$. Then we again sum up the degrees and obtain the contradiction.
\end{proof}

\begin{rem}
Unlike the approach here, the previous papers~\cite{ba1982,bf1984,pach1998,grom2010} mostly considered discrete measures concentrated on finite point sets in $\mathbb R^d$. In this case Theorems~\ref{bfbcg} and \ref{bfbcg2} hold, because we may approximate a discrete measure by an absolutely continuous measure, distributed on a set of $\delta$-balls with centers at the original concentration points. After going to the limit $\delta\to 0$ we may also assume by the standard compactness reasoning that the centers $c_\delta$ also tend to some point $c$. Then a simple argument shows that $c$ is the required point for the original discrete measure.

Direct application of the above reasoning to discrete measures is also possible, but in this case we have to deal with some accumulated error, because the number $\varepsilon$ cannot be arbitrarily small; some $k$-face of $Y$ still has to meet at least one random $(d-k)$-simplex.
\end{rem}

\begin{rem}
Imre~B\'ar\'any has noted that Theorem~\ref{bfbcg} implies the colorful Tverberg theorem\footnote[0]{Given a family of $(d+1) T(r, d)$ points in $\mathbb R^d$ colored into $d+1$ colors each containing $T(r,d)$ points, there exist $r$ disjoint ``rainbow'' $(d+1)$-tuples of points such that the corresponding $r$ convex hulls of the $(d+1)$-tuples have a common point.}~\cite{bl1992,vz1992,bmz2009} with a bad bound $T(r,d)$ of order
$$
T(r,d) \sim \frac{r}{1 - (1 - p_d)^{1/(d+1)}}\sim r (d+1)!(d+1).
$$
Of course, this bound is much worse that the known other bounds (the optimal bounds are in~\cite{bmz2009} and have order $r$), but unlike the previous known proofs this proof uses very little topology.
\end{rem}

The author thanks Arseniy~Akopyan, Imre~B\'ar\'any, J\'anos~Pach, and the two referees for the discussion and useful remarks.

\end{document}